\documentclass[12pt]{amsart}
\usepackage{amsmath, amssymb, amsthm, latexsym}
\input xypic
\newtheorem{theorem}{Theorem}
\newtheorem{proposition}[theorem]{Proposition}
\newtheorem{observation}[theorem]{Observation}

\newtheorem{lemma}[theorem]{Lemma}

\def\Proof{\medskip\noindent{\bf Proof: }}

\def\Z{\mathbb{Z}}

\def\C{\mathbb{C}}

\def\R{\mathbb{R}}
\def\C{\mathbb{C}}

\def\S{\tilde{S}}

\def\Pi{\mathbb{P}^{\infty}}

\def\qed{\hfill$\square$\medskip}

\def\Zpk{\mathbb{Z}/p^{k}}
\def\Zpk1{\mathbb{Z}/p^{k-1}}

\newcommand{\rref}[1]{(\ref{#1})}

\newcommand{\cform}[3]{\begin{array}{c}
{\scriptstyle #3}\\
#1\\
{\scriptstyle #2}\end{array}}

\newcommand{\beg}[2]{\begin{equation}\label{#1}#2\end{equation}}
\def\r{\rightarrow}

\def\S{\mathbb{S}}

\def\sl2{\widetilde{SL_{2}(\Z)}}

\title{Topological Hermitian Cobordism}
\author{Po Hu and Igor Kriz}
\thanks{The first author was supported in part by NSF grant DMS 1104348. 
The second author was supported by NSF grant DMS 1102614
}

\begin{document}

\maketitle

\begin{abstract}
Extending our method for investigating
Real cobordism (which was recently used by Hill, Hopkins and Ravenel in their
solution of the Kervaire invariant $1$ problem), we investigate the 
$RO(G)$-graded homotopy groups
of a (non-complete) $\Z/2\times \Z/2$-equivariant spectrum called
topological Hermitian cobordism. The methods of this paper may be useful in
computing the homotopy groups of other $G$-equivariant spectra where
$G\neq \Z/2$.
\end{abstract}

\section{Introduction}

A decade ago, the authors \cite{hk} studied Real cobordism $M\R$, an
$RO(\Z/2)$-graded
$\Z/2$-equivariant spectrum discovered by Landweber \cite{land}, which is
related to complex cobordism $MU$ in the same way as Atiyah's Real $K$-theory
$KR$ \cite{a} does to $K$-theory. In particular, the authors computed
the $RO(\Z/2)$-graded coefficients of $M\R$. A key step in that computation 
was the fact that $M\R$ is a {\em complete} $\Z/2$-equivariant
spectrum, which means that the spectrum of its fixed points
is equivalent to its spectrum of homotopy fixed points.
At that time, the authors thought
about the possibility of enhancing, in some way,
the $\Z/2$-action on $M\R$ or some modification of it and computing the coefficients
in hopes of obtaining more homotopy-theoretical information,
but didn't make much progress.
The subject was recently revived in a spectacular way by Hill, Hopkins and
Ravenel \cite{hhr} who solved negatively the famous Kervaire invariant $1$ problem
using an $RO(\Z/8)$-graded $\Z/8$-spectrum obtained, roughly
speaking, from
smashing four copies of $M\R$ and considering the $\Z/8$-action which combines 
a cyclic permutation of the factors with the Real action. 

\vspace{3mm}
In the present paper, we consider a more straightforward action on $M\R$ itself,
by the group $\Z/2\times\Z/2$. This action, also, gives rise to an $RO(\Z/2\times\Z/2)$-graded
spectrum (or, in the terminology of \cite{lms}, spectrum indexed over the complete
universe), which we call $M\R_{\Z/2}$. The action is defined precisely in the next section, but roughly
speaking, it combines the Real structure with a structure of ``$\Z/2$-equivariant
cobordism''. In fact, the authors were aware of this action while writing \cite{hk},
but couldn't calculate the coefficients at that time. 

\vspace{3mm}
Recently, the authors' interest in $M\R_{\Z/2}$ was renewed for a reason
unrelated to \cite{hhr}:
the construction of $M\R_{\Z/2}$ 
is a direct cobordism analogue of Karoubi's topological $\mathfrak{L}$-theory
\cite{kar},
which can also be thought of as ``$K\R_{\Z/2}$''. In a joint paper with Ormsby \cite{hko},
the authors gave a set of foundations of $G$-equivariant stable motivic
homotopy theory, and used it to give a solution to the homotopy fixed point problem
for Karoubi's algebraic Hermitian $KH$-theory for a field of
characteristic $0$, where the full force of motivic homotopy theory
can be brought to bear. \cite{hko} also contains a definition
of algebraic Hermitian cobordism $MGL\R$. Algebraic Hermitian cobordism
is still quite mysterious,
but when we specialize to the field $\R$, similarly as in \cite{kar} for
the case of $K$-theory, there is a topological version, which turns
out to coincide with $M\R_{\Z/2}$. Therefore, knowing more about $M\R_{\Z/2}$
gives information about $MGL\R$.

\vspace{3mm}
At the same time, we also realized that $M\R_{\Z/2}$ may be calculationally accessible,
and this is the subject of the present paper. Interestingly, the techniques are a little
different than we expected. In \cite{hhr}, a ``slice spectral sequence'' plays a crucial
role. This is also the case in \cite{hko}. However, in the case of \cite{hko}, the main
point was to define a motivic analogue of the ``Tate diagram'' of Greenlees and May \cite{gm},
and then use a slice spectral sequence to investigate Tate cohomology, which is more accessible
than the Borel cohomology itself.

\vspace{3mm}
In the case of the $\Z/2\times\Z/2$-spectrum $M\R_{\Z/2}$, one also has an analogue
of such ``Tate diagram'', but two things should be pointed out. First of all, the relevant
diagram is {\em not} the case for $G=\Z/2\times\Z/2$
of the Tate diagram canonically associated to $G$-equivariant spectra
for any compact Lie group by
Greenlees and May. The diagram we need is a generalized construction which brings more fully
to bear the theory of classifying spaces of families by Lewis, May and Steinberger \cite{lms},
to take a ``Tate diagram'' with respect to a $\Z/2$ subgroup of $\Z/2\times \Z/2$
(which we call $\Z/2\{h\}$), and 
then consider separately the fixed points of the $\Z/2\{h\}$-fixed point spectrum
under the quotient $\Z/2$. The other thing to realize is that $M\R_{\Z/2}$ is {\em not}
a complete spectrum, and thus our ``Tate diagram'' plays a somewhat different role than in
completion theorems. The Tate diagram was previously used in a similar way by the second author
in the simpler problem of computing the coefficients of $MU_{\Z/p}$ \cite{kriz}.

\vspace{3mm}
In the final step, we encounter $M\R^*$ cohomology groups of ($\Z/2$-fixed) stunted projective
spaces, which is also a calculation we were aware of in \cite{hk} as desirable, and couldn't do.
In the present paper, we do this calculation by taking advantage of the fact that $M\R$, again,
is a complete spectrum. We completely identify the differentials of the corresponding Borel
cohomology spectral sequence, and as a result compute an associated
graded object of the
$RO(\Z/2\times\Z/2)$-graded coefficients of $M\R_{\Z/2}$ with respect to a suitable complete
filtration (Theorem \ref{tmroff}); 
we lack a good enough ``nomenclature'' for the elements to solve all the extensions
at this point. However, for a subring graded by certain special dimensions, we do have a complete
answer as a ring (Theorem \ref{tmrtate}). 

\vspace{3mm}
The present paper is organized as follows. All of our precise statements are too technical to make
in the introduction. In the Section \ref{thc}, we establish the notation to state 
the non-calculational part of our result. We also do all the relevant equivariant stable homotopy
theory. In Section \ref{trc}, we compute the $M\R^*$ cohomology of stunted projective spaces.
In the Appendix, we say a few words on how $M\R_{\Z/2}$ relates to $MGL\R$.

\vspace{3mm}

\section{The topological hermitian cobordism spectrum and its coefficients}
\label{thc}

\noindent
{\bf Notation:}
We begin with establishing some notation. $\Z/2$-equivariant Real cobordism is a 
$\Z/2\times\Z/2$-equivariant spectrum indexed over the complete universe. This means
that we must distinguish carefully
between the different elements of $\Z/2\times \Z/2$ and its different real
irreducible representations. We will denote the non-zero elements of 
$\Z/2\times\Z/2$ by $g_\alpha$, $g_\gamma$ and $h$, and its non-trivial
real irreducible representations by $\alpha$, $\gamma$ and $\gamma\alpha$:
By definition, on $\alpha$, $h,g_\alpha$ act by minus, on $\gamma$, $h, g_\gamma$
act by minus, on $\gamma\alpha$, $g_\alpha,g_\gamma$ act by minus. We will think
of $g_\gamma$, $g_\alpha$ as ``real structures'', and of $h$ as the ``$\Z/2$-equivariant
structure''.

\vspace{3mm}
Accordingly, we will consider a complete complex universe, i.e.
$$U=(\C[\Z/2\{h\}])^\infty.$$
We denote complex conjugation by $g_\alpha$, and put $g_\gamma=hg_\alpha$. (However,
note that we can think of $g_\gamma$ as the complex conjugation.) For a $\Z/2$-equivariant
complex $\Z/2$-representation $V$, denote by $Gr(V,n)$ the space of all $n$-dimensional complex 
vector subspaces of $V$. Then as usual, there is a ``tautological'' 
$\Z/2$-equivariant complex vector bundle $\gamma^{n}$ on $Gr(V,n)$ where the fiber over an
$n$-subspace $W\subset V$ is $W$. We let
\beg{edmr1}{(M\R_{\Z/2})_V:=Gr(U\oplus V,|V|)^{\gamma_{|V|}}}
(where the superscript denotes the Thom space).
For $V\subset W\subset\subset U,$ $Gr(U\oplus V,|V|)$ is canonically embedded
into $Gr(U\oplus W,|W|)$ by adding $W-V$ (in the second summand),
and accordingly the restriction of the bundle $\gamma_{|W|}$
on $Gr(U\oplus W,|W|)$ to $Gr(U\oplus V,|V|)$ canonically splits off the equivariant
``trivial'' (i.e. induced from a point) bundle $W-V$, i.e. we get a canonical map
$$S^{W-V}\wedge (M\R_{\Z/2})_V\r (M\R_{\Z/2})_W.$$
Taking into account the whole $\Z/2\times\Z/2$-action, the resulting 
$\Z/2\times\Z/2$-equivariant spectrum indexed over the complete universe $U$
is what we denote by $M\R_{\Z/2}$. It is worth commenting that the construction
we described is obviously promoted to a $\Z/2\times\Z/2$-equivariant
symmetric spectrum, and hence $M\R_{\Z/2}$ is a $\Z/2\times \Z/2$-equivariant
$E_{\infty}$ ring spectrum indexed over the complete universe. 

\vspace{3mm}
The goal of
this paper is to compute the ``coefficients'' 
\beg{emr*}{(M\R_{\Z/2})_{k+\ell\gamma\alpha+m\alpha+n\gamma},\; k,\ell,m,n\in\Z.
}
There is one simplification which we may deduce right away. Recall \cite{cgk}
that a $\Z/2$-equivariant commutative associative
ring spectrum $E$ indexed over the complete universe
is called {\em complex-oriented} if for every $\Z/2$-equivariant finite-dimensional
complex vector bundle $\xi$ of dimension $n$ on a $\Z/2$-equivariant CW-complex $X$,
there exists a Thom class, i.e. a class
$$u_\xi\in\widetilde{E}^{2n}X^\xi$$
such that if we denote by $\theta:X^\xi\r X^\xi\wedge X_+$ the Thom diagonal,
then we have an isomorphism
\beg{eco+}{\diagram\protect\theta^*(u_\xi\otimes ?):E^{k+\ell\alpha}X
\rto^(.6)\cong & \protect\widetilde{E}^{k+\ell\alpha+2n}X^\xi\enddiagram
}
where $\alpha$ is the sign representation. (The definition really works for 
$G$-equivariant spectra for any
finite abelian group $G$, but the case of $G=\Z/2$ is the only one we need here.)

\begin{observation}
\label{oco}
When $E$ is a complex-oriented $\Z/2$-equivariant spectrum, then the coefficients
of $E$ are $2-2\alpha$-periodic. $MU_{\Z/2}$ is a complex-oriented $\Z/2$-equivariant
spectrum.
\end{observation}

\Proof
For the first statement, consider the $\Z/2$-equivariant complex bundle
$2\alpha$ over a point. For the second statement, apply classification of
equivariant complex $n$-bundles to $Gr(U\oplus n,n)$.
\qed

\vspace{3mm}
Now following Atiyah \cite{a}, a {\em $\Z/2$-equivariant Real bundle}
is a $\Z/2$-equivariant complex bundle with a $\Z/2$-equivariant antilinear
involution. Following the conventions at the beginning of this section,
we call a $\Z/2\times \Z/2$-equivariant commutative associative ring spectrum
indexed over the complete universe {\em Real-oriented} if for
every $n$-dimensional Real bundle on $\Z/2\times\Z/2$-space $X$,
there exists a Thom class, i.e. a class
$$u_\xi\in\widetilde{E}^{n(1+\gamma\alpha)}X^\xi,$$
such that we have an isomorphism
$$\diagram
\theta^{*}(u_\xi\otimes ?):
E^{k+\ell\alpha+m\gamma+r\gamma\alpha}X
\rto^{\cong} &
\widetilde{E}^{k+\ell\alpha+m\gamma+r\gamma\alpha+n(1+\gamma\alpha)}X^{\xi}.
\enddiagram
$$

\begin{lemma}
\label{lro}
When $E$ is a Real-oriented $\Z/2\times\Z/2$-equivariant spectrum,
indexed over the complete universe, the coefficients of 
$E$ are $(\gamma+\alpha-1-\gamma\alpha)$-periodic. The spectrum
$M\R_{\Z/2}$ is Real-oriented. Additionally, we have
\beg{elro1}{M\R_{k+\ell\alpha+m\gamma+n\gamma\alpha}\cong
M\R_{k+m\alpha+\ell\gamma+n\gamma\alpha}.
}
\end{lemma}

\Proof
For the first statement, consider the Real $\Z/2$-equivariant bundle 
$\alpha+\gamma$ over
a point. For the second statement, consider the $n$-dimensional
Real bundle classifying spaces 
$Gr(U\oplus n, n)$. For the last statement, note that
$\alpha$ and $\gamma$ play symmetric roles in the definition
of $M\R_{\Z/2}$.
\qed

\vspace{3mm}

Lemma \ref{lro} identifies certain dimensions in the coefficients
of $M\R_{\Z/2}$, so it reduces the set of separate ``dimensions''
we must consider. To get an further, however, we must substantially
use equivariant stable homotopy theory, as developed in
Lewis, May and Steinberger \cite{lms}.
The method we present here may well be more general, and useful
in computing the coefficients of $G$-equivariant spectra where
$G$ is a non-cyclic finite abelian group.

The strategy, in our case, is to investigate the
$\Z/2$-equivariant spectrum 
$$(M\R_{\Z/2})^{\Z/2\{h\}}.$$
This, in turn, can be computed by considering the 
``$\Z/2\{h\}$-Tate diagram'' in the language of 
Greenlees-May \cite{gm}. In our case, it is important
to note that this is {\em not} the same approach
as considering directly the $\Z/2\times \Z/2$-equivariant
Tate diagram. One must note that after taking $\Z/2\{h\}$-fixed points,
our ``Tate diagram'' retains a $\Z/2$-equivariant
action by taking 
$$\Z/2=(\Z/2\times\Z/2)/(\Z/2\{h\}).$$
Fortunately, the foundations of what we need have been completely
set up in \cite{lms}. Using the terminology of \cite{lms},
we have a $\Z/2\times\Z/2$-equivariant cofibration sequence
\beg{eeqt}{E\mathcal{F}_+\r S^0\r S^{\infty\alpha+\infty\gamma}}
where $\mathcal{F}$ is the family of subgroups disjoint with
$\Z/2\{h\}$ (recall that the {\em classifying space}
of a family $\mathcal{F}$ of subgroups of a finite group $G$
is a $G$-CW complex $E\mathcal{F}$ such that for $H\subseteq G$, $E\mathcal{F}^H$
is contractible when $H\in \mathcal{F}$ and empty otherwise). The diagram we
have in mind is
\beg{emrtate}{
\diagram
E\mathcal{F}_+\wedge M\R_{\Z/2}\rto\dto_\simeq & M\R_{\Z/2}\rto\dto &
S^{\infty\alpha+\infty\gamma}\wedge M\R_{\Z/2}\dto \\
E\mathcal{F}_+\wedge F(E\mathcal{F}_+,M\R_{\Z/2})\rto &
F(E\mathcal{F}_+,M\R_{\Z/2})\rto & 
S^{\infty\alpha+\infty\gamma}\wedge F(E\mathcal{F}_+,M\R_{\Z/2}).
\enddiagram
}
The reason the left column is an equivalence is as follows: For
a family $\mathcal{F}$ of subgroups of $G$, we have
a notion of $\mathcal{F}$-equivalence of $G$-equivariant
spectra, which means equivalence on $H$-fixed points for every
$H\in \mathcal{F}$. Now it is obvious that the natural map
$$M\R_{\Z/2}\r F(E\Z/2_+,M\R_{\Z/2})$$
is an $\mathcal{F}$-equivalence, and hence by Lemma 2.12 of 
\cite{lms}, the left hand column of \rref{emrtate} is
an $\mathcal{F}$-equivalence. However, these spectra are
homotopy equivalent to $\mathcal{F}$-CW spectra (i.e. CW-spectra
whose cells are of the form $G/H_+\wedge D^n$, $n\in \Z$, 
$H\in\mathcal{F}$), and by Theorem 2.2 of \cite{lms},
the left hand column of \rref{emrtate} is a weak equivalence).
Now applying $(?)^{\Z/2\{h\}}$ to \rref{emrtate},
we therefore obtain a (weak) homotopy pullback
\beg{emrtate3}{\diagram
(M\R_{\Z/2})^{\Z/2\{h\}}\rto\dto & 
(S^{\infty\alpha+\infty\gamma}\wedge M\R_{\Z/2})^{\Z/2\{h\}}\dto \\
F(E\mathcal{F}_+,M\R_{\Z/2})^{\Z/2\{h\}}\rto &
(S^{\infty\alpha+\infty\gamma}\wedge F(E\mathcal{F}_+, M\R_{\Z/2}))^{\Z/2\{h\}}
\enddiagram
}
It is also worth noting that by the Adams isomorphism \cite{lms} Theorem 7.1,
the homotopy fiber of the rows of \rref{emrtate3} is
$$(E\mathcal{F}_+\wedge M\R_{\Z/2})/(\Z/2\{h\}).$$

\vspace{3mm}
Now the main point of our method is that the upper left, lower left and lower
right corners of the pullback \rref{emrtate3} can be calculated
directly. In many was, in fact, the computation is analogous to \cite{kriz}; if
we do not do any suspensions by sums of copies of $\alpha$ or $\gamma$, the
complications introduced by the Real structure are in fact only minor. We
will treat this case first, in part because in this case, we have a more
precise theorem.

\vspace{3mm}
Regarding the upper right corner, the main idea is
that we can compute the $\Z/2\{h\}$-fixed point of the spectrum 
$S^{\infty\alpha+\infty\gamma}\wedge M\R_{\Z/2}$
on the prespectrum level and then take the colimit. This relies on a
result of Lewis-May-Steinberger \cite{lms} that in general, for
a normal subgroup $H$ of a (say) finite group $G$, and a based CW $G$-space
$X$, if we denote by $\mathcal{F}_H$
the family of subgroups not containing $H$, and consider the cofibration
sequence
$$(E\mathcal{F}_H)_+\r S^0\r \widetilde{E\mathcal{F}_H},$$
then we have an equivalence of $G/H$-spectra indexed over the complete universe
\beg{egeomfp}{(\widetilde{E\mathcal{F}_H}\wedge \Sigma^{\infty}_{G}X)^H
\simeq \Sigma^{\infty}_{H}X^H.
}
Note that in the case discussed in this paper, $\mathcal{F}_{\Z/2\{h\}}$ coincides with
the family of subgroups of $\Z/2\times \Z/2$ disjoint with $\Z/2\{h\}$.
Similarly as in tom Dieck \cite{td} 
(cf. \cite{kriz}), this then identifies the upper right corner
as
\beg{emrtate4}{\cform{\bigvee}{n\in\Z}{}\Sigma^{n(1+\gamma\alpha)}B\mathbb{U}_+
\wedge M\R
}
where $\mathbb{U}$ is the infinite unitary group with $\Z/2$-action
by complex conjugation. 

\vspace{3mm}

The lower left corner of \rref{emrtate3} is computed as follows: 
Consider the inclusion of universes
$$i:U^{\Z/2\{h\}}\r U.$$
Then over $U^{\Z/2\{h\}}$, we have a ``$\Z/2\{h\}$-fixed Real cobordism
spectrum with completely replicates the definition of $M\R_{\Z/2}$ with
$U$ replaced by $U^{\Z/2\{h\}}$ (so the action of $\Z/2\{h\}$ on all
spaces constituting the prespectrum is trivial). This is also essentially
the same construction as the construction of
$M\R$ (cf. \cite{hk}), so we will also denote this $U^{\Z/2\{h\}}$-indexed
$\Z/2\times\Z/2$-spectrum by $M\R$. Then we have
a standard map 
\beg{emrtate4a}{M\R\r i^*M\R_{\Z/2}
}
is an $\mathcal{F}$-equivalence, so the natural map on $\Z/2\{h\}$-fixed
points is a $\Z/2$-equivariant weak equivalence:
\beg{emrtate5}{\diagram
F(E\mathcal{F}_+,M\R)^{\Z/2\{h\}}\rto^\simeq &
F(E\mathcal{F}_+,M\R_{\Z/2})^{\Z/2\{h\}}.
\enddiagram
}
(Recall that in applying $(?)^{\Z/2\{h\}}$ to a spectrum indexed
over the complete universe, we apply $i^*$ implicitly first.)

\vspace{3mm}
The left hand side of \rref{emrtate5} is, by definition,
\beg{eclass}{F(E\mathcal{F}/(\Z/2\{h\})_+,M\R).}
However, we have 
$$E\mathcal{F}/(\Z/2\{h\})\simeq B_{\Z/2}(\Z/2)$$
(the right hand side is, by definition, the classifying space of $\Z/2$-equivariant
principal $\Z/2$-bundles). In the next section, we will show that the natural
inclusion
$$B_{\Z/2}(\Z/2)\r B\S^1$$
(where $\S^1$ denotes the unit sphere in $\C$ with $\Z/2$-action by complex
conjugation) induces in $M\R$-cohomology (in dimensions $k+\ell\gamma\alpha$,
$k,\ell\in\Z$) the map
\beg{emrtate6}{M\R^*B\S^1=M\R^*[[u]]\r M\R^*[[u]]/[2]_F(u)
}
where $F$ is the universal formal group law, and $u$ is a variable in
dimension $-1-\gamma\alpha$ (just as in \cite{hk}, in an effort to prevent
constant confusions, all gradings are homological).

\vspace{3mm}
Combining \rref{emrtate4}, \rref{emrtate5}, \rref{emrtate6}, we get that the
diagram of the
coefficients of the upper right, lower left and lower right corner of the
diagram \rref{emrtate3} in dimensions $k+\ell\gamma\alpha$, $k,\ell\in \Z$
is
\beg{emrtate7}{\diagram
&M\R^*[u,u^{-1}][b_1,b_2,...]\dto\\
M\R^*[[u]]/[2]_F(u)\rto & (M\R^*[[u]]/[2]_F(u))[u^{-1}]
\enddiagram
}
The generators $b_k$ are in dimensions $k(1+\gamma\alpha)$ and just as in
\cite{kriz}, Theorem 1.1, the images of $b_k$ are the coefficients
of $x^k$ in $x+_F u$. Also, just as in \cite{kriz}, the images of the
upper right and lower left corners span the lower right corner, and we get

\vspace{3mm}

\begin{theorem}
\label{tmrtate}
The bigraded module of coefficients
$$((M\R_{\Z/2})_{k+\ell\gamma\alpha})_{k,\ell\in\Z}$$
is isomorphic to the pullback of the diagram of rings \rref{emrtate7}.
\end{theorem}
\qed

\vspace{3mm}

Of course, we are not done with computing the $RO(\Z/2\times\Z/2)$-graded
coefficients of $M\R_{\Z/2}$, since we must also discuss suspensions by
linear combinations, with coefficients in $\Z$, of $\alpha$ and $\gamma$.
Let us realize first precisely which coefficients are left to compute. We claim,
in fact, that it suffices to compute
\beg{emrr}{(M\R_{\Z/2})_{k+\ell\gamma\alpha+m\alpha}, \; m< 0.}
To see this, consider an arbitrary coefficient \rref{emr*}. By Lemma \ref{lro},
we may assume, say, that $m\leq n$. Using the periodicity in Lemma \ref{lro},
however, we may then add $n$ to $k,\ell$ and subtract $n$ from $m,n$ and the 
resulting isomorphic coefficient is either of the form \rref{emrr} (when $m<n$),
or of the form covered in Theorem \ref{tmrtate}(when $m=n$).

\vspace{3mm}
Now the idea is to consider \rref{emrr} one $m$ at a time. This means replacing
the diagram \rref{emrtate3} by
\beg{emrtate3a}{\diagram
(\Sigma^{m\alpha}M\R_{\Z/2})^{\Z/2\{h\}}\rto\dto & 
(S^{\infty\alpha+\infty\gamma}\wedge M\R_{\Z/2})^{\Z/2\{h\}}\dto \\
F(\Sigma^{-m\alpha}E\mathcal{F}_+,M\R_{\Z/2})^{\Z/2\{h\}}\rto &
(S^{\infty\alpha+\infty\gamma}\wedge F(E\mathcal{F}_+, M\R_{\Z/2}))^{\Z/2\{h\}},
\enddiagram
}
$-m>0$. (We realize that the right column is periodic with respect to suspending
by a multiple of $\alpha$.) Therefore, the main task remaining is to compute
\beg{emrtate4aa}{F(\Sigma^{-m\alpha}E\mathcal{F}_+,M\R_{\Z/2})^{\Z/2\{h\}}.}
Using the ``splitting'' \rref{emrtate4a}, we get that \rref{emrtate4aa}
is equivalent to 
\beg{emrtate4ab}{F(\Sigma^{-m\alpha}E\mathcal{F}_+,M\R)^{\Z/2\{h\}}
\simeq F((\Sigma^{-m\alpha}E\mathcal{F}_+)/(\Z/2\{h\}),M\R_{\Z/2}).
}
Note that the right hand side is simply a $\Z/2$-spectrum. But we can do better.
Recall \rref{eclass}. It follows that 
\beg{emrrtate1}{(\Sigma^{-m\alpha}E\mathcal{F}_+)/(\Z/2\{h\})\simeq 
(B_{\Z/2}(\Z/2))^{-m\gamma_1}
}
where $\gamma^1$ is the ``tautological'' $1$-dimensional real bundle
on $B_{\Z/2}(\Z/2)$, induced by the sign representation of $\Z/2$ on $\R$.
Additionally, again, if we consider the {\em fixed} $\Z/2$-space
$\R P^\infty$, we have a canonical inclusion
$$(\R P^\infty)^{-m\gamma_1}\r (B_{\Z/2}(\Z/2))^{-m\gamma_1},$$
which, since $M\R$ is a complete spectrum, induces an equivalence
\beg{emrrtate2}{F((B_{\Z/2}(\Z/2))^{-m\gamma_1},M\R)\simeq 
F((\R P^\infty)^{-m\gamma_1},M\R).
}
Note that one commonly refers to $(\R P^\infty)^{-m\gamma_1}$ as
the {\em stunted projective space} $\R P^{\infty}_{-m},$ so
our problem is reduced to computing 
\beg{emrrtate3}{\widetilde{M\R}^{*}(\R P^{\infty}_{-m}),\; -m\geq 0}
(where $\R P^{\infty}_{-m}$ is considered as a fixed $\Z/2$-space),
and the bottom row of diagram \rref{emrtate3a}, which, in view of 
naturality and Theorem \ref{tmrtate}, follows from knowing
\beg{emrrtate4}{\widetilde{M\R}^{*}(\R P^{\infty}_{-m})\r
M\R^{*}(\R P^{\infty}),
}
induced by the ``$0$-section'' map 
$$\R P^\infty\r \R P^{\infty}_{-m}, \; m>0.$$
Unfortunately, we cannot determine \rref{emrrtate4} completely. The
main problem is that we have no obvious ``pretty'' expression
for \rref{emrrtate3} akin to \rref{emrtate6}. What we do have
is a complete computation of the Borel cohomology spectral sequence (BCSS) for
\rref{emrrtate4}, and the map to the BCSS for $M\R^*\R P^\infty$.
Ordinarily, one would think this is not sufficient information
to determine the coefficients of the upper left corner of \rref{emrtate3}
or even its associated graded object, since taking associated graded objects
does not preserve pullbacks. However, in the present case, the pullbacks
we are dealing with are rather special, and therefore it is possible 
to state a theorem on this level of precision. Let us begin with a
preliminary lemma.

\vspace{3mm}
Let $F^i,\;i\geq 0$ be a decreasing filtration on an abelian
group $X$. We will call this filtration {\em complete}
if the natural map
\beg{excompl}{X\r \cform{\lim}{\leftarrow}{} X/F^i X
}
is an isomorphism.

\vspace{3mm}

\begin{lemma}
\label{lpb}
Suppose we have a cartesian co-cartesian square of abelian groups
\beg{elpb0}{
\diagram
A\rto\dto & B\dto^f\\
C\rto_g &D.
\enddiagram
}
(Recall that this is the same thing as a cartesian square, or pullback diagram,
such that $Im(f)+Im(g)=D$.) Suppose further we have complete decreasing filtrations 
$F^i,\;i\geq 0$ on $B,C,D$ such that 
\beg{elpba}{f(F^i B)+g(F^i C)= F^i D.
}
Then forming the pullback
\beg{elpbb}{
\diagram
F^iA\rto\dto & F^iB\dto^{F^if}\\
F^iC\rto_{F^ig} &F^i{D},
\enddiagram
}
$(F^iA)_i$ is a complete decreasing filtration on $A$, and we also
have a pullback of the corresponding associated graded objects
\beg{elpbc}{
\diagram
E_0A\rto\dto & E_0B\dto^{E_0f}\\
E_0C\rto_{E_0g} &E_0{D}.
\enddiagram
}
\end{lemma}

\Proof
We may view as a cartesian co-cartesian square \rref{elpb0} as 
a short exact sequence
\beg{elpb1}{\diagram 0\rto & A\rto & B\oplus C\rto^(.6){f\oplus g} & D\rto & 0.\enddiagram
}
Assuming \rref{elpbb}, we have a short exact sequence embedded into \rref{elpb1}
\beg{elpb2}{\diagram 0\rto & F^iA\rto & F^iB\oplus F^iC\rto^(.65){F^if\oplus F^ig} 
& F^iD\rto & 0.\enddiagram
}
Therefore, we have a quotient short exact sequence
\beg{elpb3}{\diagram 0\rto & A/F^iA\rrto && B/F^iB\oplus C/F^iC\rrto^(.65){f/F^if\oplus g/F^ig} 
&& D/F^iD\rto & 0.\enddiagram
}
Now interpreting \rref{elpb3} as a pullback again, taking 
inverse limits over $i$ and using 
the commutation of categorical limits gives the completeness of the filtration $(F^i A)_i$.
Taking the quotient of \rref{elpb2} by the same short exact sequence with $i$ replaced
by $i+1$ gives the last statement.
\qed

\vspace{3mm}
Consider now the fixed point inclusion of $\Z/2\times\Z/2$-spaces
\beg{ebdef}{b:S^0\r S^{\gamma\alpha}.
}

\vspace{3mm}
\begin{theorem}
\label{tmroff}
There is a complete filtration on \rref{emrr} such that the associated graded object
is isomorphic to the pullback of the diagram
\beg{emrtate7a}{\diagram
&E_0M\R^*[u,u^{-1}][b_1,b_2,...]\dto\\
E_0\widetilde{M\R}^*(\R P^{\infty}_{-m})\rto & E_0(M\R^*[[u]]/[2]_F(u))[u^{-1}]
\enddiagram
}
where in the upper right, lower left
and lower right corner of \rref{emrtate7a},
$E_0$ denotes the associated graded object of the decreasing filtration by powers
of the ideal $(b)$. The lower horizontal arrow of \rref{emrtate7a} is as
computed in Theorems \ref{trpodd}, \ref{trpeven} below.
\end{theorem}

\vspace{3mm}
\noindent
{\bf Comment:} It is {\em not} being asserted that the filtration of \rref{emr*}
mentioned in the Theorem is the filtration by powers of the ideal $(b)$.

\vspace{3mm}
\Proof
(with the exception of the last statement, which is proved in the next section).
The strategy is to show that the $RO(\Z/2)$-graded coefficients
of diagram \rref{emrtate3a}
satisfy the hypotheses of Lemma \ref{lpb}. To this end, note that 
the lower left corner is $M\R^*$-cohomology of a spectrum (reduced
cohomology of a space, actually), so the filtration by powers
of $(b)$ is simply the BCSS filtration, and the statement follows
from the convergence of the corresponding BCSS (to be completely
precise, the BCSS is only conditionally convergent, so the proof
follows from our complete computation of the differentials, which
is done in the next section). 

To prove completeness of the filtration by powers of $(b)$ in the lower right
and upper right corner of the diagram \rref{emrtate7a},
our strategy is to show that any infinite series of elements in the same
dimension divisible by an increasing power of $b$ can be written as an
infinite series of elements divisible by an increasing power of $b$ using 
only $u^m,\; m\geq -n$ with constant $n$. Such series then converge in
$u^{-n}M\R^*[[u]]/[2]_F(u)$, which is a shift of $M\R^*\R P^\infty$, where
the series converges by the completeness of $M\R$ (since the filtration by
powers of $(b)$ is the filtration associated with the BCSS).

The statement asserted in the last paragraph follows, in effect, from
dimensional considerations. Recall from \cite{hk} (or verify directly)
that the cokernel of the canonical inclusion 
\beg{emrcok}{\Z[a]/(2a)\r M\R_{k+\ell\alpha}}
is $0$ when $k+\ell<0$. Since $u^{-1}$ is of dimension $1+\alpha$ (from the
point of view of $M\R$, the dimension is $1+\gamma\alpha$ from the point of
view of our entire calculation), if we have a homogeneous 
infinite series as mentioned in the last
paragraph involving powers of $u$ not bounded below, then
the $M\R_*$-coefficients of those powers must be in the image of \rref{emrcok}. 
(In the case of the upper right corner, recall that $b_k$ is in dimension
$k(1+\alpha)$.)
One easily sees however that as we increase the power of $a$ in \rref{emrcok},
adding no multiple of $1+\alpha$ can bring the elements into the same dimension.
This proves the required completeness statement. 

To complete verifying the assumptions of Lemma \ref{lpb}, it remains
to verify \rref{elpba}. This is done as follows: Denote by $\gamma_{1}^{\C}$ the
$\Z/2$-equivariant Real line bundle on $B_{\Z/2}(\Z/2)$ given by letting
the generator of $\Z/2$ act on $\C$ by $-1$. Then we have the ``$0$-section map''
\beg{ebbound}{\iota:(B_{\Z/2}(\Z/2))^{\gamma_1}\r (B_{\Z/2}(\Z/2))^{\gamma_{1}^{\C}}.
}
Applying $M\R^*$, and inverting $u$, we get a diagram
\beg{ebbound1}{
\diagram
u^{m}M\R^*[[u]]/[2]_F(u)\rto^{\iota^*}\drto_{\kappa} &
\widetilde{M\R}^*\R P^{\infty}_{-m}\dto\\
& (M\R^*[[u]]/[2]_F(u))[u^{-1}].
\enddiagram
}
In the upper right corner, we use the equivalence \rref{emrrtate2} again,
and in the upper left corner we use the Real orientation of $-m\gamma_{1}^{\C}$.
In fact, this also implies that $\kappa$ is the ordinary localization map,
so its image together with the image of the vertical map \rref{emrtate7}
span $(M\R^*[[u]]/[2]_F(u))[u^{-1}]$. This implies \rref{elpba} for $i=0$.
Note, however, that since $F^i=b^iF^0$ in all the three terms of
\rref{emrtate7a}, the general case follows.
\qed

\vspace{3mm}

\section{The Real cobordism of stunted projective spaces}
\label{trc}

The following result is well known, but we restate it to make the exposition
self-contained:

\begin{lemma}
\label{lfil}
Let $E^{*}_{**}$ be a spectral sequence (graded homologically). Suppose
we have numbers $p,s$ and a morphism of spectral sequences
$$\phi:E^{\prime*}_{**}\r E^{*}_{**}$$
which is an isomorphism on $E^{s}_{m,q}$-terms for $m\leq p$, and $E^{s}_{m,q}=0$ for
$m>p$. Then for $r\geq s$, $m\leq p$ 
\beg{elf1}{E^{\prime r}_{*,*}\cong E^{r}_{*,*}\oplus \cform{\bigoplus}{s\leq i<r}{} Im (d^{i}_{c})}
where $d^{i}_{c}$ is the restriction of $d^i$ to 
$$E^{i}_{>p,*}\r E^{i}_{\leq p,*}.$$
Furthermore, in the isomorphism \rref{elf1}, $\phi$ corresponds to the projection to the first
summand.
Similarly, suppose 
we have a morphism of spectral sequences
$$\psi:E^{*}_{**}\r E^{*}_{\prime**}$$
which is an isomorphism on $E^{s}_{m,q}$-terms for $m> p$, and $E^{s}_{m,q}=0$ for
$m\leq p$. Then for $r\geq s$, $m\geq p$ 
\beg{elf2}{E^{\prime r}_{*,*}\cong E^{r}_{*,*}\oplus \cform{\bigoplus}{s\leq i<r}{} Coim (d^{i}_{c}).}
Furthermore, in the isomorphism \rref{elf2}, $\psi$ corresponds to the injection of the first
summand.
\end{lemma}

\Proof
Induction on $r$. For $r=s$, the statement is assumed. Consider now, say, $\phi$.
Assuming the induction hypothesis for a given
$r$, the second summand by definition may be non-zero only in filtration degreess $>p-r+1$, and hence
$\phi$ is an isomorphism on $E^{r}_{\leq p-r+1,*}$-terms. This means that $E^{\prime r}_{\leq p-r,*}$,
which can be taken as the 
codomain of $d^r$-differentials
in $E^{\prime r}$, 
is mapped by $\phi$ monorphically, and hence $\phi d_r=d_r\phi$ determines $d_r$. This
implies the statement with $r$ replaced by $r+1$, hence completing the induction. The
proof of the statement for $\psi$ is analogous.
\qed

\vspace{3mm}
For example, this lemma (applied with $s=1$) implies that quite generally, 
the differentials of a Borel homology and Borel cohomology spectral
sequence with respect to a cyclic group are precisely the restrictions
of the differentials in the corresponding Tate spectral sequence
(see also \cite{gm}).

\vspace{3mm}
Let $\R P^\infty$ denote the infinite real projective space with trivial $\Z/2$-action.
observe that 
\beg{erpfixed}{\R P^\infty=(B\S^1)^{\Z/2}
}
where $\S^1$ is the unit sphere in $\C$ with the $\Z/2$-action of complex conjugation.
Recall \cite{hk} that by the theory of Real-oriented spectra,
$$BP\R^*B\S^1=BP\R^*[[u]],$$
where the dimension of $u$, graded homologically, is $-1-\alpha$.

\begin{proposition}
\label{prp}
The relation \rref{erpfixed} induces
an isomorphism
$$BP\R^{*}\R P^\infty\cong BP\R^*[[u]]/[2]_F(u).$$
\end{proposition}

\Proof
If we denote by $\S^\infty$ the unit sphere in $\C^\infty$ with $\Z/2$-action
by complex conjugation, then $\S^\infty/(\Z/2)$ is the unit sphere of the
Real line bundle $(\gamma_1)^2$ on $B\S^1$ where $\gamma^1$ is the 
tautological bundle. Thus, we have a cofibration sequence
$$\S^\infty/(\Z/2)\r B\S^1\r (B\S^1)^{(\gamma^1)^2},$$
which gives an isomorphism
$$BP\R^*(\S^\infty/(\Z/2))\cong BP\R^*[[u]]/[2]_F(u).$$
On the other hand, notice that
\rref{erpfixed} factors through the natural equivariant embedding
\beg{erpinc}{\R P^\infty =(\S^\infty)^{\Z/2}/(\Z/2)\subset \S^\infty/(\Z/2)}
which is an equivalence non-equivariantly and hence induces an isomorphism
in $BP\R$-cohomology, since $BP\R$ is a complete spectrum.
This completes the proof.
\qed

Note that \rref{erpinc} is not an equivalence equivariantly: in fact, the
right hand side is $B_{\Z/2}(\Z/2)$, the classifying space for $\Z/2$-equivariant
$\Z/2$-principal bundles.

\vspace{3mm}
Even though \ref{prp} gives a complete calculation of $BP\R^*(\R P^\infty)$,
for the purposes of comparison with stunted projective spaces,
we need to study its Borel cohomology spectral sequence. Let us begin
by studying the $[2]_F$-series in $BP\R_*$ further. Even though
one may think this is the same as the $[2]_F$-series in $BP_*$, there are some
fundamental differences. For example, the natural map
$$BP^*[[u]]/(\frac{[2]_F(u)}{u})\r BP^*[[u]]/(\frac{[2]_F(u)}{u})[u^{-1}]$$
is an inclusion. This is not the case when we replace $BP^*$ by $BP\R^*$.
If we choose a natural number $n$, then in $BP\R^*[[u]]$,
there exists a natural number $C$ such that 
\beg{eprp1}{a^{2^n-1}([2]_F(u)+Cv_nu^{2^{n+1}-1})\cong a^{2^n-1}
(v_nu^{2^n}+v_{n+1}u^{2^{n+1}})\mod u^{2^{n+
1}-1}.
}
We denote the left hand side by $a^{2^n-1}\phi_nu^{2^n}$. Note that 
therefore $\phi_nu^{2^n}$  makes good sense, although $\phi_n$
does not (at least not in $BP\R^*[[u]]$).
We shall also put 
$$w_na^{2^n-1}u^{2^n}:=a^{2^n-1}(\phi_n-v_n).$$
Again, $w_n a^{2^n-1}$ makes good sense, while $w_n$ does not.

\vspace{3mm}
\begin{lemma}
\label{lprp}
The Borel cohomology spectral sequence for $BP\R^*(\R P^\infty)$
has $E^1$-term 
$$BP^*[a][\sigma,\sigma^{-1}][[u]]/[2]_F(u)$$
where, to conform with \cite{hk}, everything is graded homologically,
$|v_n|=(2^n-1)(1+\alpha)$, $|a|=-\alpha$, $|\sigma|=\alpha-1$,
$|u|=-1-\alpha$. The differentials are
\beg{edif1}{d (x\sigma^{-2^n})=xa^{2^{n+1}-1}v_n,}
where
\beg{edif1a}{x\in\Z/2[v_n,v_{n+1},...][a][\sigma^{\pm 2^{n+1}}]u^j
}
with
$$\;0\leq j<2^n,$$
(when $n=0$, $\Z/2[v_0,v_1...]$ is replaced by $BP_*$) and
\beg{edif2}{d(y\sigma^{-2^n})=ya^{2^{n+1}-1}w_n,}
where
\beg{edif2a}{y\in\Z/2[v_{n+1},...][a][\sigma^{\pm 2^{n+1}}]u^j 
}
with
$$\;j\geq 2^n.$$
\end{lemma}

\Proof
The statement about the $E^1$-term is obvious. The fact that $u$ is a permanent
cycle follows from Proposition \ref{prp}. The differential \rref{edif1} follows
from the fact that the Borel cohomology spectral sequence (BCSS) for $BP\R^*\R P^\infty$
is a module over the BCSS for $BP\R^*$. The differential \rref{edif2} is
a rewriting or \rref{edif1} by
the discussion in the paragraph preceeding the Lemma. There can be no further
differentials by Lemma \ref{lfil}
since, as one easily checks, after introducing the differentials
\rref{edif1}, \rref{edif2},
the corresponding ``Tate spectral sequence'' (i.e. the spectral
sequence obtained by inverting $a$)
converges to its correct target, $H\Z/2^*\R P^\infty[a,a^{-1}]$.
\qed

\begin{lemma}
\label{lbot}
When $q=2p$, and 
\beg{elbot1}{\ell=2^{n_1}-2^{n_2}+....-2^{n_{q}},}
$n_1>...>n_{q}$ or 
$q=2p+1$ and 
\beg{elbot2}{\ell=2^{n_1}-2^{n_2}+...+2^{n_{q}}-1,}
$n_1>...>n_q$, then in the BCSS for $BP\R^*\R P^{\infty}_{2\ell+1}$,
\beg{elbot3}{du^{\ell+1}=\sigma^{2^{n_q}}v_qa^{2^{n_q+1}-1}u^{\ell+1}.
}
\end{lemma}

\Proof
As an induction hypothesis, we
claim that the kernel of the Tate cohomology spectral sequence (TCSS) map 
induced by the projection
$$\R P^\infty\r \R P^{\infty}_{2\ell+1}$$
where $\ell \cong 2^k \mod 2^{k+1}$ or $\ell \cong 2^k+1\mod 2^{k+1}$
on the $2^m-1$'st page, $m\leq k$, is a free 
$\Z/2[v_m,v_{m+1},...][a,a^{-1}][\sigma^{\pm 2^{m}}]$-module
$M_m$ on the generators 
\beg{elbott1}{\phi_m u^{\ell+i},\; 1\leq i\leq 2^m.
}
For $m=1$, this is a statement about the $E^1$-term. Suppose it is true
for a given $m$. Then first of all, $d^{<2^{m+1}-1}$ of $\sigma^{-2^n}u^{\ell+1}$
must be contained in $M_m$ (because the differential is $0$ in the TCSS associated
to $BP\R^*\R P^\infty$). However, we see that $d^1$ is excluded by
explicit formula, and for dimensional
reasons, $d^{>1}$ can only arise
from a shift of power of $\sigma$, so $d^{2^{m+1}-1}$ is the lowest possible
differential (corresponding to shift by $\sigma^{2^m}$), so $d^{<2^{m+1}-1}$
cannot occur.

Now from considering the map of TCSS associated with
\beg{elprdrp}{\Sigma^{2\ell+1}M\Z/2\r \R P^{\infty}_{2\ell+1},}
we see that for $m<\ell$, 
\beg{elprdif}{d^{2^{m+1}-1}\sigma^{-2^m}u^{\ell+1}=v_ma^{2^{m+1}-1}u^{\ell+1} 
+\cform{\sum}{1<i\leq 2^m}{}
x_i u^{\ell+i}
a^{2^{m+1}-1}}
where $x_i\in \phi_m\Z/2[v_n,v_{n+1},...]$
(the power of $\sigma$ is excluded for dimensional reasons). 
But then by the module structure over the TCSS for $BP\R^*$,
$$d^{2^{m+1}-1}u^{\ell+1}=\cform{\sum}{1<i\leq 2^m}{}\sigma^{2^m}x_i u^{\ell+i}
a^{2^{m+1}-1}.$$
However, for the lowest $i$ for which $x_i\neq 0$, we then have
by \rref{elprdif} 
$$d^{2^{m+1}-1}(\sigma^{2^m}x_i u^{\ell+i}
a^{2^{m+1}-1})\neq 0 \mod u^{\ell+i+1},$$
which is a contradiction with $d^{2^{m+1}-1}$ being a differential. 
Therefore, $x_i=0$ for all $i$ and
$$d^{2^{m+1}-1}\sigma^{-2^m}u^{\ell+1}=v_m
a^{2^{m+1}-1}u^{\ell+1}$$
and
$$d^{2^{m+1}-1}u^{\ell+1}=0.$$
This, in turn, implies the induction step. 

Now let us consider $d^{2^{k+1}-1}$. Then by \rref{elprdrp},
$$d^{2^{\ell+1}-1}u^{\ell+1}=\phi_\ell\sigma^{2^\ell} a^{2^{\ell+1}-1} 
+\cform{\sum}{1<i\leq 2^\ell}{}
x_i u^{\ell+i}
a^{2^{m+1}-1},$$
$$d^{2^{\ell+1}-1}(u^{\ell+1}\sigma^{-2^\ell})=w_\ell a^{2^{\ell+1}-1} 
+\cform{\sum}{1<i\leq 2^\ell}{}
x_i u^{\ell+i}\sigma^{-2^\ell}
a^{2^{\ell+1}-1}.$$
Again, however, for the lowest $i$ for which $x_i\neq 0$, we conclude
that then 
$$d^{2^{\ell+1}-1}(x_i u^{\ell+i}\sigma^{-2^\ell}
a^{2^{\ell+1}-1})\neq 0 \mod u^{\ell+i+1},$$
which is a contradiction with $d^{2^{\ell+1}-1}$ being a differential.
This implies that $x_i=0$ for all $i$, and the statement of the Lemma.
\qed

\vspace{3mm}
Now let $\ell$ be as in \rref{elbot1} or \rref{elbot2}.
Let
$$\epsilon(\ell):=2\ell-2^{n_1}+1.$$
Let further
$$\ell_1:=\ell,$$
$$\ell_{i+1}:= 2^{n_i}-\ell_i-1,\; i<q.$$
Set 
$$\epsilon_i:=1+\ell+\epsilon(\ell_q)+...+\epsilon(\ell_{i+1}).$$

\begin{lemma}
\label{lbot1}
For $0\leq i\leq q,$ in the BCSS for $BP\R^*\R P^{\infty}_{2\ell+1}$
we have
$$d(u^{\epsilon_i})=\sigma^{2^{n_i}}v_{n_i}a^{2^{n_i+1}-1}u^{\epsilon_i}.$$
\end{lemma}

\Proof
One notes that if one puts
$$\ell^\prime:=2^{n_1}-2^{n_2}+...-2^{n_i} \;\text{if $i$ is even}$$
and
$$\ell^\prime:=2^{n_1}-...+2^{n_i}+1 \;\text{if $i$ is odd},$$
then 
$$\ell^\prime\leq \epsilon(\ell_q)+...+\epsilon(\ell_{i+1})<\ell^\prime +2^{n_i}-1.$$
If $i$ is odd, then $\ell^\prime\geq \ell$, so we have a map
$$\R P^{\infty}_{2\ell +1}\r \R P^{\infty}_{2\ell^\prime +1},$$
and we may therefore deduce the differential from Lemma \ref{lbot}
applied to $\ell^\prime$ in place of $\ell$.

When $i$ is even, then $\ell^\prime\leq \ell$, so we have a map in
the opposite direction
$$\R P^{\infty}_{2\ell^\prime +1}\r \R P^{\infty}_{2\ell +1},$$
so the TCSS for $BP\R^*\R P^{\infty}_{2\ell^\prime +1}$ {\em detects}
the differential, but a priori, we don't know that there isn't an
error term. However, in this case, we may proceed by reversed induction
on $i$: if the statement is correct for $i^\prime >i$, 
then we may use the differentials already proved and the module
structure of the TCSS $E^\prime$ corresponding to $BP\R^*\R P^{\infty}_{2\ell^\prime +1}$,
and the TCSS $E$ corresponding to $BP\R^*\R P^{\infty}_{2\ell +1}$,
and prove that the map
$$E^{2^{n_i+1}-1}\r E^{\prime 2^{n_i+1}-1}$$
is injective, and hence the differential in the target determines the differential
in the source (see the statement of Theorem \ref{trpodd} 
which spells out this calculation explicitly).
\qed

\vspace{3mm}
\begin{theorem}
\label{trpodd}
The BCSS for $BP\R^*\R P^{\infty}_{2k+1}$ has $E^1$-term
$$u^{k+1}BP_*[a][\sigma,\sigma^{-1}][[u]]/(\frac{[2]_F(u)}{u}),$$
($|u|=-1-\alpha$). The differentials are as follows: For each $i\geq k+1$,
there are numbers $s_{k,i}\leq t_{k,i}\in \{1,2,...,\infty\}$ such 
that we have \rref{edif1} with \rref{edif1a} for
$$0\leq p<s_{k,i}$$
or 
$$p=t_{k,i},$$
and
\rref{edif2} with \rref{edif2a} for
$$s_{k,i}\leq p<t_{k,i}.$$
Additionally, we have the differential
\beg{edif3}{dz=v_{p}a^{2^{p+1}-1}\sigma^{2^p},\;p=t_{k,i}
}
with
\beg{edif3a}{z\in u^i\Z/2[v_p,v_{p+1},...][a][\sigma^{\pm 2^{p+1}}].
}
These are the only differentials in this BCSS.

The numbers $s_{k,i}, t_{k,i}$ are determined as follows:
For $k=2^{\ell}-1$, put 
$$t_{k,i}=\ell \;\text{for 
$k+1\leq i\leq 2k,$}$$ 
$$t_{k,i}=\infty \;\text{for $i>2k$}.$$
We put 
$$s_{k,i}=j\; \text{for $k+2^{j-1}\leq i<k+2^j,$ 
$1\leq j\leq \ell$},$$ 
and 
$$s_{k,i}=j\; \text{for $i>2k$, 
$2^j\leq i<2^{j+1}.$}$$

For $k$ not of the form $2^\ell-1$, let $k+1\leq 2^\ell
\leq 2k$. Then
$$
\begin{array}{l}
s_{k,i}=s_{2^\ell-k-1,2^\ell-2k-1+i}\\
t_{k,i}=t_{2^\ell-k-1,2^\ell-2k-1+i}
\end{array}
\;
\text{for $i<2^\ell$,}
$$
and if we let
$$2^\ell-k\leq 2^m<2^{\ell+1}-2k-2,$$
then for $2^\ell\leq i\leq 2k,$
put
$$t_{k,i}=\ell, \; s_{k,i}=j$$
for
$$m<j\leq \ell,\; 2k+1-2^\ell+2^{j-1}\leq i<2k+1-2^\ell+2^j.$$
For $i>2k,$
put 
$$t_{k,i}=\infty,\; s_{k,i}=j$$
for $2^j\leq i<2^{j+1}$.
\end{theorem}

\Proof
The differentials \rref{edif3} with \rref{edif3a}
follows from Lemma \ref{lbot1}. The remaining differentials follow
from the module structure over the BCSS for $BP \R^* \R P^\infty$. 
Considering the corresponding TCSS, one observes that these differentials
get the correct answer for $\widetilde{BP\R}^*\R P^{\infty}_{2k+1}=
H\Z/2^*\R P^{\infty}_{2k+1}$, and hence no other differentials can
occur.
\qed

\vspace{3mm}

The result for $BP\R^*\R P^{\infty}_{2k}$ is analogous, although the
pattern is a little different.

\vspace{3mm}
\begin{theorem}
\label{trpeven}
The BCSS for $BP\R^*\R P^{\infty}_{2k}$ has $E^1$-term
$$u^{k}BP_*[a][\sigma,\sigma^{-1}][[u]]/([2]_F(u)),$$
($|u|=-1-\alpha$). The differentials are as follows: For each $i\geq k$,
there are numbers $a_{k,i}\leq b_{k,i}\in \{0,2,...,\infty\}$ such 
that we have \rref{edif1} with \rref{edif1a}, \rref{edif2} with
\rref{edif2a} and \rref{edif3} with \rref{edif3a} occur
as in Theorem \rref{trpodd} with $s_{k,i}$ replaced by $a_{k,i}$
and $t_{k,i}$ replaced by $b_{k,i}$.

The numbers $a_{k,i}, b_{k,i}$ are determined as follows:
For $k=2^{\ell}$, put 
$$b_{k,i}=\ell \;\text{for 
$k\leq i< 2k,$}$$ 
$$b_{k,i}=\infty \;\text{for $i\geq 2k$}.$$
We put 
$$a_{k,i}=j\;\;\;\; \parbox{2in}{for $i=k$ and
$0=j$, or $k+2^{j-1}\leq i<k+2^j,$ 
$1\leq j\leq \ell$},$$ 
and 
$$a_{k,i}=j\;\; \text{for $i\geq 2k$, 
$2^j\leq i<2^{j+1}.$}$$

For $k$ not of the form $2^\ell$, let $k\leq 2^\ell
< 2k$. Then
$$
\begin{array}{l}
a_{k,i}=a_{2^\ell-k,2^\ell-2k+i}\\
b_{k,i}=t_{2^\ell-k,2^\ell-2k+i}
\end{array}
\;
\text{for $i<2^\ell$,}
$$
and if we let
$$2^\ell-k\leq 2^m<2^{\ell+1}-2k,$$
then for $2^\ell\leq i< 2k,$
put
$$b_{k,i}=\ell, \; a_{k,i}=j$$
for
$$m<j\leq \ell,\; 2k-2^\ell+2^{j-1}\leq i<2k-2^\ell+2^j.$$
For $i>2k,$
put 
$$b_{k,i}=\infty,\; a_{k,i}=j$$
for $2^j\leq i<2^{j+1}$.
\end{theorem}

\vspace{3mm}
The proof is completely analogous to the proof of Theorem \ref{trpodd},
so we only record the main steps. Lemma \ref{lbot} is replaced by

\vspace{3mm}
\begin{lemma}
\label{lebot}
When
\beg{elebot1}{\ell=2^{n_1}-2^{n_2}+...\pm 2^{n_q},
}
(the sign is $+$ resp. $-$ when $q$ is odd resp. even),
$n_1>...>n_{q-1}>n_q+1$, then
in the BCSS for $BP\R^*\R P^{\infty}_{2\ell}$,
there is a differential
\beg{elebot3}{d u^\ell=u^\ell\sigma^{2^{n_q}}v_{n_q}a^{2^{n_q}-1}.
}
\end{lemma}

\Proof
Analogous to the proof of Lemma \ref{lbot}. The principal difference
is that \rref{elprdrp} is replaced by the inclusion of the bottom cell
\beg{elpredrp}{S^{2\ell}\r \R P^{\infty}_{2\ell}.
}
\qed

\vspace{3mm}
To get the analogue of Lemma \ref{lbot1}, let $\ell$ be as in \rref{elebot1},
and let 
$$\delta(\ell):=2\ell=2^{n_1},$$
and let
$$\ell_1:=\ell,$$
$$\ell_{i+1}:=2^{n_i}-\ell_i,\; i<q.$$
Set 
$$\delta_i:=\ell+\delta(\ell_{i+1})+...+\delta(\ell_q).$$
The analogue of Lemma \ref{lbot1} then reads

\vspace{3mm}
\begin{lemma}
\label{lebot1}
For $0\leq i\leq q,$
$$d(u^{\delta_i})=
\sigma^{2^{n_i}}v_{n_i}a^{2^{n_i+1}-1}u^{\delta_i}.$$
\end{lemma}

The proof is analogous to the proof of Lemma \ref{lbot1}, and the proof
of Theorem \ref{trpeven} is analogous to the proof of Theorem \ref{trpodd}.

\vspace{3mm}

\noindent
{\bf Example:} Let us compute the Borel cohomology spectral sequence for
$BP\R^*\R P^{\infty}_{101}$. This is the case of Theorem \ref{trpodd} with
$k=50$. For the purposes of the theorem, we write
$$k=2^6-2^4+2^2-2^1.$$
We refer to the differentials \rref{edif1} as {\em vertical}, the differentials \rref{edif2}
as {\em horizontal} and the differentials \rref{edif3} as {\em special}. 

In this terminology, we get special differentials
$$d_3\sigma^{0\mod 4}u^{51},$$
$$d_7\sigma^{0\mod 8} u^{52},$$
$$d_{31}\sigma^{0\mod 32}u^i,\; 53\leq i\leq 63,$$
$$d_{127}\sigma^{0\mod 128} u^i,\; 64\leq i\leq 100.$$
In addition, we have vertical differentials
$$d_{15}\sigma^{8\mod 16},\; 53\leq i\leq 56,$$
$$d_{63}\sigma^{32 \mod 64}u^i,\; 64\leq i\leq 68$$
and horizontal differentials
$$d_1\sigma^{1\mod 2}, \; i\geq 51,$$
$$d_3\sigma^{2\mod 4},\; i\geq 51,$$
$$d_7 \sigma^{4\mod 8}u^i,\; i\geq 52,$$
$$d_{15}\sigma^{8\mod 16},\; i\geq 56,$$
$$d_{31}\sigma^{16\mod 32},\; i\geq 53,$$
$$d_{63}\sigma^{32\mod 64},\; i\geq 69,$$
$$d_{127}\sigma^{64\mod 128},\; i\geq 64.$$
Differentials $d\sigma^ju^i$ for $i\geq 101$
are the same as in the BCSS for $BP\R^* \R P^\infty$.

\vspace{3mm}

\noindent
{\bf Example:} Let us compute the Borel cohomology spectral sequence for
$BP\R^*\R P^{\infty}_{100}$. This is the case of Theorem \ref{trpeven} with
$k=50$. For the purposes of the theorem, we write
$$k=2^6-2^4+2^1.$$

In the above terminology, we get special differentials
$$d_3\sigma^{0\mod 4}u^{i},\; i=50,51,$$
$$d_{31}\sigma^{0\mod 32}u^i,\; 52\leq i\leq 63,$$
$$d_{127}\sigma^{0\mod 128} u^i,\; 64\leq i\leq 99.$$
In addition, we have vertical differentials
$$d_1\sigma^{1\mod 2}u^{50},$$
$$d_{15}\sigma^{8\mod 16},\; 52\leq i\leq 55,$$
$$d_{63}\sigma^{32 \mod 64}u^i,\; 64\leq i\leq 67$$
and horizontal differentials
$$d_1\sigma^{1\mod 2},\; i\geq 51,$$
$$d_3\sigma^{2\mod 4},\; i\geq 50,$$
$$d_7 \sigma^{4\mod 8}u^i,\; i\geq 52,$$
$$d_{15}\sigma^{8\mod 16},\; i\geq 56,$$
$$d_{31}\sigma^{16\mod 32},\; i\geq 52,$$
$$d_{63}\sigma^{32\mod 64},\; i\geq 68,$$
$$d_{127}\sigma^{64\mod 128},\; i\geq 64.$$
Differentials $d\sigma^j u^i$ for $i\geq 100$
are the same as in the BCSS for $BP\R^* \R P^\infty$.

\section{Appendix: The topological realization of algebraic hermitian cobordism}
\label{sapp}

Let us start with a warm-up case, Karoubi's topological Hermitian $K$-theory 
\cite{kar} $\mathfrak{L}$, which is the topological realization of his
algebraic Hermitian $K$-theory over the field $\R$. 
When considered as a $\Z/2\times \Z/2$-equivariant 
spectrum indexed over the complete universe, this is represented as follows:
One chooses a complex universe $U_B\cong \C^\infty$ 
together with a Real structure $\overline{?}$
and a hyperbolic Real form $B$. Then the $\mathfrak{L}^0$ is represented by
the $\Z/2\times \Z/2$-equivariant space $BGL(U_B)$, the classifying
space of the group $GL(U_B)$ which is the direct limit of groups of the
form $GL(V_B)$ where $V_B$ are finite-dimensional Real vector spaces
endowed with the $\Z/2\times \Z/2$-action where one generator acts by
complex conjugation, and the other by 
\beg{econjinv}{A\mapsto (A^T)^{-1}} 
where $(?)^T$
denotes the adjoint matrix with respect to the symmetric bilinear form $B$:
$$B(Ax,y)=B(x,A^Ty).$$
The link with $\Z/2$-equivariant Real $K$-theory is obtained as follows.
Since we assumed we have a Real structure on $V_B$, we have
a unique splitting
$$V_B=(V_B)^+\oplus (V_B)^-$$
such that $(V_B)^+$, $(V_B)^-$ are $B$-orthogonal, and $B$ is positive
definite on $(V_B)^+$ and negative definite on $(V_B)^-$. Then define
a positive definite symmetric bilinear form $B_+$ to be equal
to $B$ on
$(V_B)^+$ and to $-B$ on $(V_B)^-$.
Now consider the Hermitian form on $V_B$ defined by
\beg{ehermb}{\langle u, v\rangle= B_+(u,\overline{V}).}
Let $U(V_B)$ be the group of unitary matrices with respect to this Hermitian
form. 
Then we have a
map
\beg{eunit}{U(V_B)\r GL(V_B)
}
which is easily seen to be a weak $\Z/2\times\Z/2$-equivalence. Moreover, define
a $\Z/2$-action on $V_B$ such that the generator $h$ acts by $1$ on $(V_B)^+$ and
by $-1$ on $(V_B)^-$. The key point then is that restricted to
$U(V_B)$, \rref{econjinv} is complex conjugation, conjugated by $h$.
Passing to direct limit, this gives an equivalence of the $0$-space of
$\Z/2$-equivariant Real $K$-theory to $\mathfrak{L}$.
One checks that the periodicity maps given by Karoubi \cite{kar} coincide,
up to homotopy, with the $K\R_{\Z/2}$-periodicity maps.

\vspace{3mm}
To treat topological Hermitian cobordism, which is the $\Z/2\times \Z/2$-equivariant
spectrum indexed over the complete universe which is the 
topological realization of the $\Z/2$-equivariant motivic spectrum
$MGL\R$ over the field $\R$, we must enrich \rref{eunit} to
``Thom spaces''. To this end, recall \cite{hko} that the ``Thom space''
in that case is obtained as the homotopy cofiber of the map from the
``unit sphere bundle'' to the base space. The reason of the quotation
marks is that the we cannot take the complement of the $0$-section in the vector bundle,
since this would not have a $\Z/2$-action. Instead, the ``unit sphere'' is 
emulated by the quadric
\beg{equad}{B(x,y)=1, \; x,y\in V_B
}
which has both a $GL_(V_B)$ action where $A$ acts by $A$ on the $x$-coordinates and by
$(A^T)^{-1}$ on the $y$ coordinates. Note that projection to the $x$-coordinates (or the
$y$-coordinates) is a motivic homotopy equivalence onto $V_B-\{0\}$.

For the purposes or comparison of the topological realization of this with 
$\Z/2$-equivariant Real cobordism, we take the topological subspace $S_+$ of 
\rref{equad} consisting of all points of \rref{equad} where 
\beg{equad1}{y=\overline{h(x)}.}
By projecting to the $x$-coordinates, this is homeomorphic to the set of all points
of $V_B$ satisfying
\beg{equad2}{B_+(x,\overline{x})=1,}
which is the unit sphere. Thus, the inclusion of $S_+$ to the quadric \rref{equad}
is a (non-equivariant homotopy equivalence). Now one sees, however, that 
the restriction of the $GL(V_B)$-action on \rref{equad} restricts
to an action of $U(V_B)$ on $S^+$, and that, in effect, the $\Z/2\ltimes U(V_B)$-space
$S_+$ is isomorphic to the unit sphere of $V_B$ with respect to the positive-definite
Hermitian form
\rref{ehermb}. On the other hand, the inclusion of $S^+$ to \rref{equad}
is a $\Z/2\ltimes U(V_B)$-homotopy equivalence.
Thus, we have exhibited homotopy equivalences of the homotopy cofiber of the projection
$$B(S_+,U(V_B),*)_+\r BU(V_B)_+$$
both to the $V_B$-space of the topological realization of $MGL\R$ over $\R$ and
of $M\R_{\Z/2}$. One easily checks that the structure maps coincide as well.

\end{document}